\documentclass{elsart3-1}

\usepackage{amssymb}

\usepackage[english,francais]{babel}

\newtheorem{e-proposition}[theorem]{Proposition}

\newtheorem{e-definition}[theorem]{Definition\rm}

\def\P{{\mathbb P}}
\def\C{{\mathbb C}}

\def\rond{\kern 1pt{\scriptstyle\circ}\kern 1pt}

\def\sdir_#1^#2{\mathrel{\mathop{\kern0pt\oplus}\limits_{#1}^{#2}}}
\input amssym.def
\input xy
\xyoption{all}
\def\Z{{\Bbb Z}}
\def\C{{\Bbb C}}
\def\P{{\Bbb P}}
\def\Z{{\Bbb Z}}
\def\pp{{\Bbb P}^2}
\def\cn{\mathop{\rm Bir}(\pp)}

\def\og{\leavevmode\raise.3ex\hbox{$\scriptscriptstyle\langle\!\langle$~}}
\def\fg{\leavevmode\raise.3ex\hbox{~$\!\scriptscriptstyle\,\rangle\!\rangle$}}

\begin{document}

\begin{frontmatter}

\selectlanguage{english}
\title{On Cremona transformations of prime
order}

\vspace{-2.6cm}

\selectlanguage{francais}
\title{Sur les transformations de Cremona d'ordre premier}
\selectlanguage{english}
\author{Arnaud Beauville}
\ead{beauville@math.unice.fr} 
\address{Laboratoire J.-A. Dieudonn\'e, UMR 6621 du CNRS,
Universit\'e de Nice, Parc Valrose, F-06108 Nice cedex 2}
\author{J\'er\'emy Blanc}
\ead{Jeremy.Blanc@math.unige.ch}
\address{Section de Math\'ematiques, Universit\'e de Gen\`eve,  2-4 rue du
Li\`evre, Case postale 240, CH-1211 Gen\`eve 24}
\begin{abstract}\par
We prove that  an automorphism of order 5 of the Del Pezzo surface
$S_5$ of degree 5 is conjugate through a birational map
$S\dasharrow \pp$ to a linear automorphism of $\pp$.
This completes the classification of conjugacy classes of elements
of prime order in the Cremona group.

\vskip 0.5\baselineskip 

\selectlanguage{francais}

\noindent{\bf R\'esum\'e}
\vskip 0.5\baselineskip\par
Nous d\'emontrons qu'un automorphisme d'ordre 5 de la surface de Del
Pezzo $S$ de degr\'e 5 est conjugu\'e via une application birationnelle 
$S\dasharrow \pp$ \`a un automorphisme lin\'eaire de  $\pp$.
Cette observation  compl\`ete la classification des classes de conjugaison
d'\'elements d'ordre premier dans le groupe de Cremona.

\end{abstract}
\end{frontmatter}
\selectlanguage{english}
\section{The result}
\par The {\em Cremona group} $\cn $ is the group of
birational transformations of $\pp_{\C}$, or 
equivalently the group of $\C$-automorphisms of the field
$\C(x,y)$. There is 
an extensive classical literature about this group, in particular
about its finite subgroups; we refer to the introduction
of \cite{dF} for a list of references. 
\par The classification of conjugacy classes of elements
of prime order $p$ in $\cn$ has been given a modern treatment
in \cite{B-B} for $p=2$ and in \cite{dF} for $p\ge 3$. Let us recall the main
results. The linear transformations of given order are contained
in a unique conjugacy class (see the Proposition below). Apart from these
there are three
 families of conjugacy classes of involutions (the
famous de Jonqui\`eres, Bertini and Geiser involutions), then two
 families of conjugacy classes of order 3, given by
automorphisms of special Del Pezzo surfaces of degree 3 and 1
respectively\footnote{A birational transformation $\sigma $ of a rational
surface $S$ defines in a natural way a conjugacy class in $\cn$, 
namely the set of elements $\pi \rond\sigma\rond \pi ^{-1}   $ where
$\pi $ runs through the set of birational maps
$S\dasharrow\pp$.}. Every transformation of prime order
$\ge 7$ is conjugate  to a linear automorphism. As for conjugacy
classes of order 5, there is at least one  family, given by automorphisms
of a special Del Pezzo surface of degree 1 (example E3 in \cite{dF}).  
\par The following result  completes
the classification:\par\smallskip 
\noindent{\bf Theorem}$.-$ {\it Every conjugacy class of order
$5$ in $\cn$ is either of the above type, or is the class containing the
  linear automorphisms.}\smallskip 
\par Most of the work has been done already in \cite{dF}, where it is
proved  that there can be only one more 
conjugacy class of order 5, namely that of
an automorphism $\sigma $ (of order 5)
of the Del Pezzo surface $S_5\subset\P^5$.
We will prove in the next section that $\sigma $ is
conjugate to a linear automorphism. 
\par The theorem, together with the classification  in \cite{dF}, has
the following consequence:
\smallskip \par
\noindent{\bf Corollary}$.-$ {\it A birational transformation of prime
order  is not conjugate to a linear automorphism if and only if it fixes
some  non-rational curve}.\smallskip 
\par Observe that this does not hold for transformations of composite
order: the automorphism $\sigma $ of the  cubic surface
$x^3+y^3+z^3+t^3=0$ in $\P^3$ given by $$\sigma (x,y,z,t)=(y,x,z,\rho
t)\quad
\hbox{with }\ \rho ^3=1,\ \rho \not=1$$has only 4 fixed points, 
while $\sigma ^2$ fixes the elliptic curve $t=0$. Thus the
conjugacy class of $\sigma$ in $\cn$ cannot contain a linear
automorphism.

\section{The proof}
 Let us first recall why linear automorphisms of the same
(finite) order are conjugate:
\par\smallskip
\noindent{\bf Proposition}$.-$ {\it Let $n$ be a positive integer. Two
linear automorphisms of order $n$ are conjugate in $\cn$}.
\vskip5pt \par
\noindent\textit{Proof} : Let $T=\C^*\times \C^*$ be the standard
maximal torus of $PGL_3$. We can view $T$ as a Zariski open subset of
$\pp$, so its automorphism group $GL(2,\Z)$ embeds
naturally in $\cn$:  an element $\pmatrix{a&b\cr c&d}$ of  
$GL(2,\Z)$ corresponds to the Cremona transformation
$(x,y)\mapsto (x^ay^b,x^cy^d)$. The image of $GL(2,\Z)$ 
 in $\cn$ normalizes $T$, and its action on $T$ by conjugation is
just the original action of $GL(2,\Z)$ on $T$. 
\par We want to prove that two elements of order
$n$ of $T$ are conjugate under this action. The kernel of the
multiplication by $n$ in $T$ is $(\Z/n)^2$; any element of
order $n$ in the $\Z/n$-module $(\Z/n)^2$ is part of a basis,
which implies that it is conjugate to $(1,0)$ under the action of
$SL(2,\Z/n)$. Since the natural map $SL(2,\Z)\rightarrow
SL(2,\Z/n)$ is surjective, this proves our result.\qed\par\medskip

\noindent\textit{Proof of the theorem} : 
\par Let $V$ be a 6-dimensional vector space, and $S\subset\P(V)$ a
Del Pezzo surface of degree 5. It admits 
an automorphism $\sigma $  of order 5, unique
up to conjugation, which comes from a linear automorphism
$s$ of $V$. 
\par The automorphism $\sigma $ has at least one fixed point
$p$ (otherwise we would have a smooth \'etale covering 
$S\rightarrow S/\langle
\sigma \rangle$ of degree 5, so that
$1=\chi(\mathcal{O}_{S})=5\chi(\mathcal{O}_{S/\langle
\sigma \rangle})$, a contradiction).  The tangent plane
$\P (T_p)$   to $S$ at $p$ corresponds to a  3-dimensional vector
subspace
$T_p\subset V$. Since $s$ preserves $T_p$, it 
induces a linear automorphism  $\bar \sigma $
 of $\pp=\P(V/T_p)$; we have a
commutative diagram
$$\xymatrix  @M=6pt{S \ar[r]^{\sigma
}\ar@{-->}[d]_{\pi }&  S\ar@{-->}[d]^{\pi }\\
\pp\ar[r]_{\bar\sigma}& \pp}$$
where $\pi :S\dasharrow \pp$ is the projection from
$\P(T_p)$.
\par We claim that $\pi $ is birational, so that it conjugates
$\sigma$ to the linear automorphism $\bar \sigma $ of $\pp$. 
Identify $S$ with $\pp$ blown up at 4 points
$p_1,\ldots ,p_4$ in general position.  The map
$\pi :S\dasharrow
\pp$ is given by the linear system $|\mathcal{C}|$ of cubics passing
through
$p_1,\ldots ,p_4$ and
 singular at  $p$. Note that the 5 points $p_1,\ldots ,p_4,p$ are
again in general position (that is, no 3 of them are collinear).
Indeed this means  that $p$ does not belong to any line in
$S$. But
$S$ contains 10 lines, on which
$\langle \sigma \rangle$ acts with 2 orbits; if $p$ lies on a line
it must lie on $5$ of them, which is impossible. 
\par Thus $\pi $ is the \textit{de Jonqui\`eres transformation} of order 3, a
classical instance of birational transformation of $\pp$ (see e.g.\
\cite{S-R}, VII, 7.2). This achieves the
proof of the theorem.\qed\bigskip

\noindent\textit{Remark}$.-$ We can make the result completely explicit.
We use the coordinates system on $\pp$ such that
$(p_1,\ldots ,p_4)$ is a projective frame; as in \cite{dF}, we can view
$\sigma$ as the lift over 
$S$ of the birational transformation  $(x:y:z)\mapsto (x(z-y):
z(x-y): xz)$ of $\pp$. This transformation has
two fixed points $(\omega :1:\omega ^2)$, with
$\omega ={1\over 2}(1\pm\sqrt 5)$. Choosing  one of these fixed
points we find that
$\sigma$  is conjugate by the involutive transformation 
$$(x:y:z)\mapsto \bigl((x-\omega y)(y-z)(z-\omega x):(\omega
^{-1} x- y)(\omega ^2y-z)(z-x):(x-y)(\omega ^2y-z)(z-\omega x)
\bigr)$$
to the linear automorphism of $\pp$ given by the matrix
$$\pmatrix{0&1&-\omega ^{-2}\cr -\omega ^{-1} &1&0\cr 0&1&0 }\ .$$


\begin{thebibliography}{00}

\bibitem[B-B]{B-B} L. Bayle, A. Beauville: {\sl Birational involutions of
$\pp$}. Kodaira's issue, Asian J.\
Math.\ {\bf 4} (2000),  11--17.

\bibitem[dF]{dF} T. de Fernex: {\sl On planar 
Cremona maps of prime order}. Nagoya Math.\ J.\ {\bf 174}
(2004). Preprint {\tt math.AG/0302175}.

 \bibitem[S-R]{S-R} J. Semple, L. Roth: {\sl Introduction to Algebraic
Geometry}.\ Oxford, at the Clarendon Press, 1949.

\end{thebibliography}
\end{document}